\documentclass[10pt]{amsart}

\usepackage{enumerate}
\usepackage{latexsym}
\usepackage[pdftex]{graphicx}
\usepackage{amssymb}
\usepackage{xcolor}
\usepackage{stackrel}
\usepackage[normalem]{ulem}

\newtheorem{lemma}{Lemma}

\newcommand{\bl}{\begin{lemma}}
\newcommand{\el}{\end{lemma}}

\newcommand{\bp}{\begin{proof}}
\newcommand{\ep}{\end{proof}}

\begin{document}

\title[Centers of mass of convex bodies]{On the maximal distance between the centers of mass of a planar convex body and its boundary}
\markright{Centers of mass of convex bodies}

\author[F~ Nazarov]{Fedor Nazarov}
\address{Department of Mathematical Sciences, Kent State University,
	Kent, OH 44242, USA} \email{nazarov@math.kent.edu}

\author[D.~Ryabogin]{Dmitry Ryabogin}
\address{Department of Mathematical Sciences, Kent State University,
	Kent, OH 44242, USA} \email{ryabogin@math.kent.edu}

\author[V.~Yaskin]{Vladyslav Yaskin}
\address{Department of Mathematical and Statistical Sciences, University of Alberta,
	Edmonton, Canada} \email{vladyaskin@math.ualberta.ca}

\maketitle

\begin{abstract}
We prove that the length of the projection of the vector joining the centers of mass of  a convex body on the plane and of its boundary to an arbitrary direction  does not exceed $\frac{1}{6}$ of the body width in this direction. It follows that the distance between these centers of mass does not exceed  $\frac16$ of the diameter of the body and  $\frac{1}{12}$ of its boundary length. None of those constants can be improved.
\end{abstract}

\section{Introduction}

Let $\Omega$ be a compact convex set with non-empty interior in ${\mathbb R^d}$, $d\ge 2$. How far apart can the centers of mass $c(\Omega)$ and $c(\partial\Omega)$ of the body $\Omega$ and its boundary $\partial\Omega$ be relative to some natural linear size of $\Omega$?  

This question was posed in \cite[A25, page 36]{CFG} and still remains mostly unresolved. In this article we consider the case when $d=2$ and the ``size" is the width $w_{\Omega}(\theta)$ of $\Omega$ in some fixed direction $\theta$. We show that
$$
\Bigl|\langle c(\partial\Omega) -  c(\Omega),\theta\rangle\Bigr|\le \frac{1}{6}w_{\Omega}(\theta).
$$
It immediately implies that
$$
\Bigl| c(\partial\Omega) -  c(\Omega)   \Bigr|\le \frac{1}{6}d({\Omega})\le \frac{1}{12}\mbox{length}(\partial \Omega),
$$
where $d({\Omega})$ and $\mbox{length}(\partial \Omega)$ are the diameter and the perimeter of $\Omega$ respectively.
The constants in these estimates cannot be improved.

The proof is a mixture of elementary analytic and geometric tools and we believe that, despite being somewhat technical in places, it still contains several ideas and tricks worth presenting to the reader.

\section{Example}

\begin{figure}[ht]
  \centering
	\includegraphics[height=4.5in]{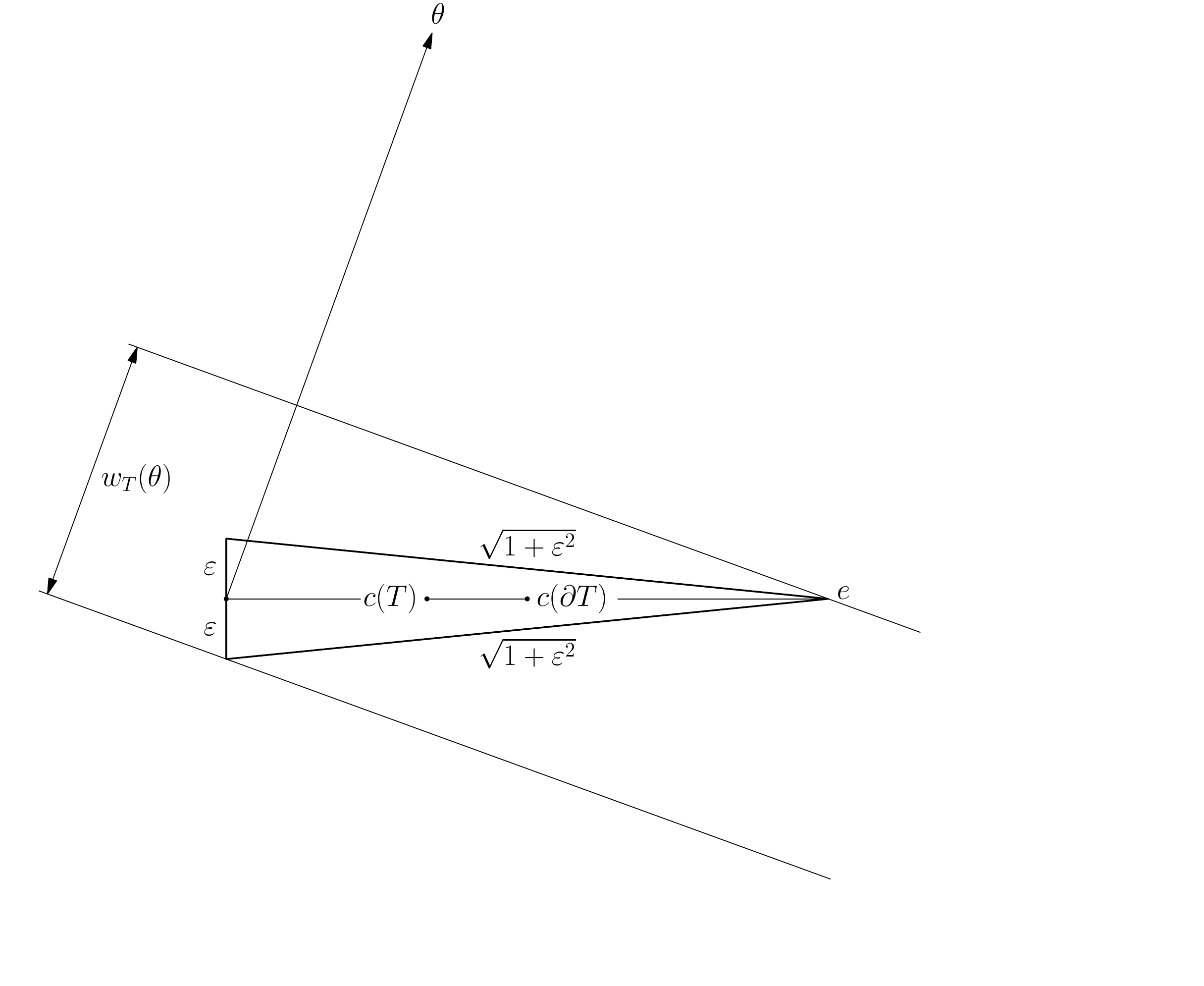}
	\caption{ }
	\label{fig1}
\end{figure}

Let $\varepsilon>0$ be a very small number. Let $T$ be the isosceles triangle with sides $\sqrt{1+\varepsilon^2}$, $\sqrt{1+\varepsilon^2}$,  and $2\varepsilon$  (see Figure~\ref{fig1}).
Let $e$ be the unit vector parallel to the vector $c(\partial T)-c(T)$. 
The distance  from the center of mass $c(T)$ of $T$ to its base  of length $2\varepsilon$ is $\frac{1}{3} $. On the other hand, the distance from the center of mass of $\partial T$ to the same line is   
$$
\frac{\sqrt{1+\varepsilon^2}}{2\varepsilon +2\sqrt{1+\varepsilon^2}}=\frac12(1+\varepsilon^2- \varepsilon \sqrt{1+\varepsilon^2}).
$$
Thus, 
$$\langle c(\partial T)-c(T),\theta\rangle = \left( \frac16+\frac12 \varepsilon^2-  \frac{\varepsilon}{2} \sqrt{1+\varepsilon^2}\right)\langle e,\theta\rangle.
$$
Also observe that the width of $T$ in every direction $\theta$ with $\langle e,\theta\rangle$ not too close to $0$ (as shown on Figure \ref{fig1}) equals 
$$w_T(\theta) = |\langle e,\theta\rangle| + \varepsilon \sqrt{1-\langle e,\theta\rangle ^2}.$$
Therefore,
$$|\langle c(\partial T)- c(T),\theta\rangle| =  \frac{ (\frac16+\frac12 \varepsilon^2-  \frac{\varepsilon}{2} \sqrt{1+\varepsilon^2})  |\langle e,\theta\rangle|}{ |\langle e,\theta\rangle| + \varepsilon \sqrt{1-\langle e,\theta\rangle ^2}} w_T(\theta),$$
where the coefficient on the right-hand side can be made as close to $\frac{1}{6}$ as one wishes.

Note now that 
$$
|c(\partial T)-c(T)|=  \frac16+\frac12 \varepsilon^2-  \frac{\varepsilon}{2} \sqrt{1+\varepsilon^2},
$$
while  we have $d(T)=\sqrt{1+\varepsilon^2}$ and  $\,\textrm{length}(\partial\Omega)=2\sqrt{1+\varepsilon^2}+2\varepsilon$, so the ratios
$\dfrac{|c(\partial T)-c(T)|}{d(T)}$ and $\dfrac{|c(\partial T)-c(T)|}{\textrm{length}(\partial T)}
$
can be made arbitrarily close to $\frac{1}{6}$ and $\frac{1}{12}$ respectively as $\varepsilon \to 0$.

\section{Outline of the proof of the upper bound}
We shall assume without loss of generality that $\theta = e_1$ is the first coordinate vector. Let $\Omega$ be any polygon without vertical sides (the general case can be obtained from this one by standard approximation arguments). We shall choose the coordinates on the $x$-axis so that the length $\ell(t)$ of the  cross-section $E_t=\{(x,y)\in \Omega: x=t\}$ is maximal at $x=0$ and the projection of $\Omega$ to the real axis is $[-1,\omega]$ for some $\omega>0$. 

Put $\Omega_t=\Omega\cap \{(x,y):x\le t\},$ for $t\in[-1,\omega]$, so $\Omega_{-1}$ is a single point and $\Omega_\omega=\Omega$. 
Let $A(t)$ and $P(t)$ be the area and the perimeter of $\Omega_t$ respectively (note that the right vertical side $E_t$ of $\Omega_t$ is included into $P(t)$). 
Let $c_a(t)$ and $c_p(t)$ be the $x$-coordinates of $c(\Omega_t)$ and $c(\partial \Omega_t)$ respectively.

\begin{figure}[ht]
	\centering
	\includegraphics[height=4in]{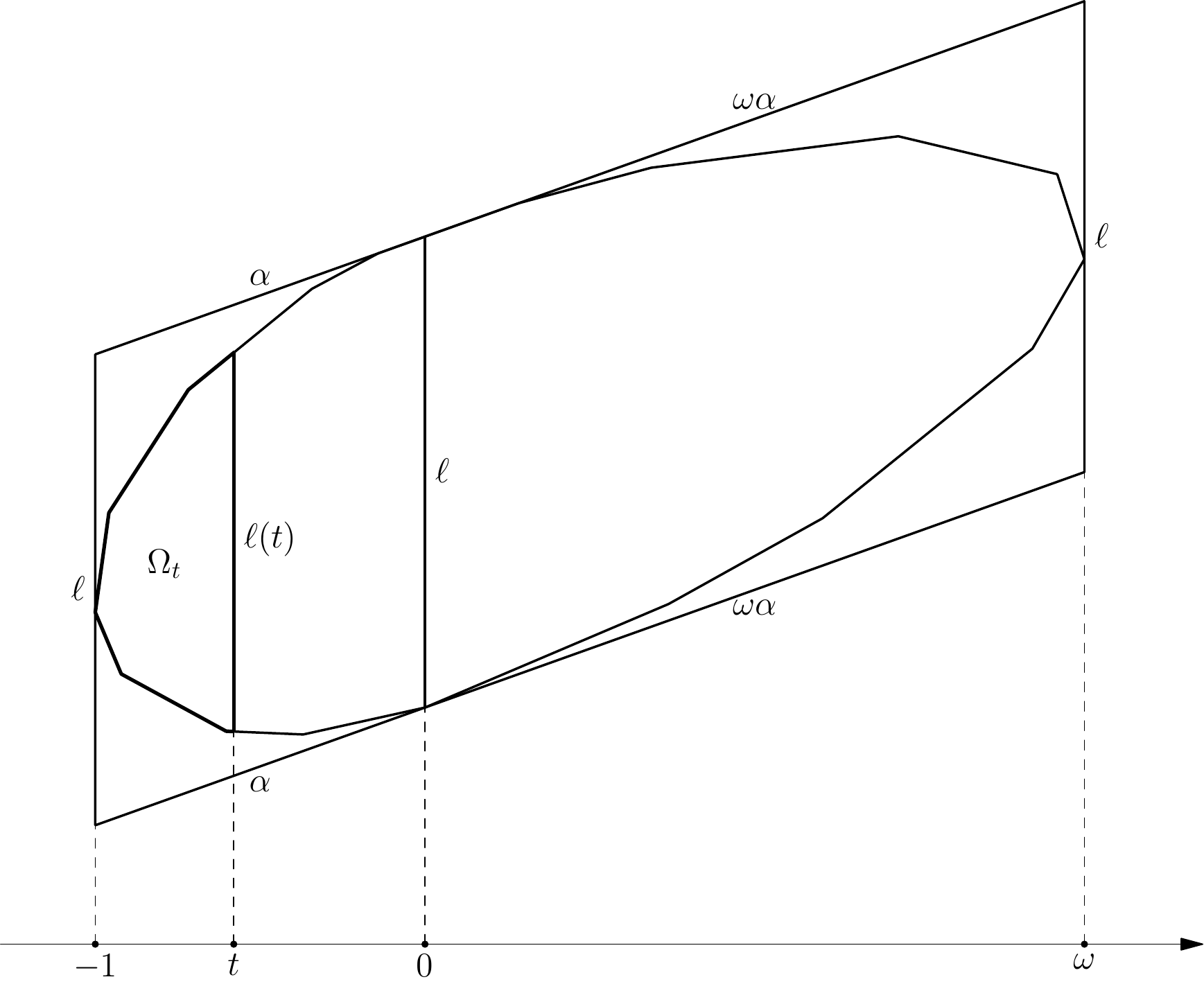}
	\caption{ }
	\label{fig2}
\end{figure}

Since $\ell:=\ell(0)=\max_t\ell(t)$, there are support  lines $L_-$ and $L_+$  to $\partial \Omega$ at the endpoints of the cross-section $E_0$ that are parallel, so $\Omega$ is inscribed into the parallelogram bounded by these two support lines and the vertical lines  $x=-1$ and $x=\omega$. The vertical sides of this parallelogram have length $\ell$. Let $\alpha(1+\omega)$ ($\alpha\ge 1$) be the common length of the two ``slanted sides" (see Figure \ref{fig2}).

The proof will consist of two main (and largely independent) parts. In Part 1, we shall establish the inequality 
\begin{equation}\label{**}
c_p(0)-c_a(0)\le \frac16.
\end{equation}
In Part 2, we will derive the final inequality 
\begin{equation}\label{***}
c_p(\omega)-c_a(\omega)\le \frac{1+\omega}6
\end{equation}
from (\ref{**}). 

Note that (\ref{**}) by itself is not sharp in the sense that the constant $\frac16$ on the right hand side can certainly be improved a bit and we suspect that the best value is actually $\frac{1}{12}$. However, this does not affect the sharpness of the final inequality because it is asymptotically  attained as $\omega\to +\infty$. 

Before these main parts, we will prove several auxiliary lemmas. We suspect that most of them are well-known but we could not locate them in the literature.

\section{Six elementary lemmas}

\bl\label{dfedsm22}
If $A(0)=B\ell$, then $c_a(0)\in \left[-\frac12, -\frac{B}{2}\right].$
\el

\begin{figure}[ht]
	\centering
	\includegraphics[height=4in, clip=true, viewport=0in 0in 5in 6in]{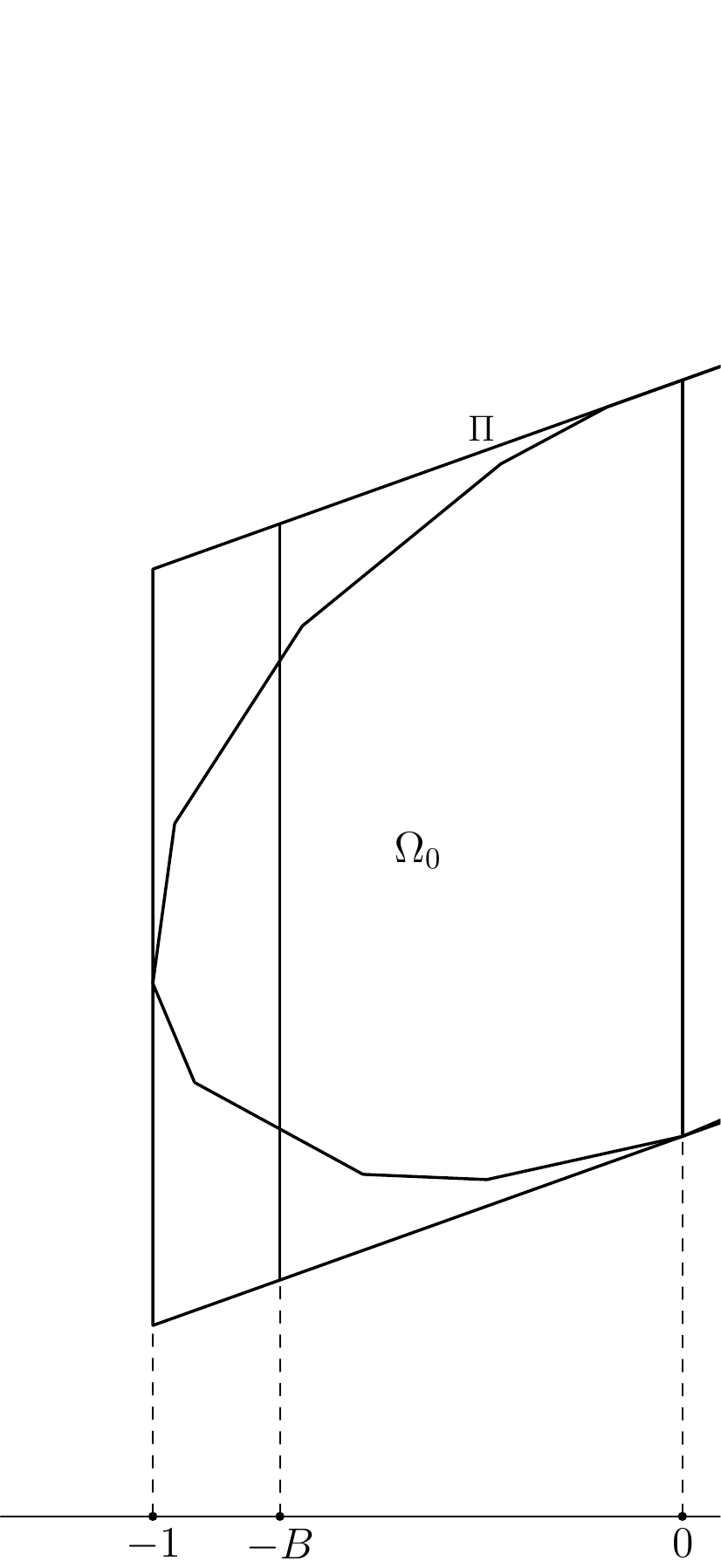}
	\caption{ }
	\label{fig3}
\end{figure}

\bp Let $\Pi$ be the parallelogram bounded by $L_+$, $L_-$ and the vertical lines  $x=-B$ and $x=0$ (see Figure~\ref{fig3}).

Then $\mathrm{Area}(\Pi)=A(0)$, so  $\mathrm{Area}(\Pi\setminus \Omega_0)=\mathrm{Area}(\Omega_0\setminus \Pi)$. However, the set $\Pi\setminus \Omega_0$ is lying to the right of the vertical line  $x=-B$ and the set $\Omega_0\setminus \Pi$ is lying to the left of the same line. Thus,
\begin{equation}\label{1.37}\langle c(\Omega_0\setminus \Pi), e_1\rangle \le -B\le \langle c(\Pi\setminus \Omega_0), e_1\rangle.
\end{equation}
We also have 
$$c(\Omega_0) = \frac{1}{A(0)}\left[\mathrm{Area}(\Omega_0\setminus \Pi)\cdot c(\Omega_0\setminus \Pi) + \mathrm{Area}(\Omega_0\cap \Pi)\cdot c(\Omega_0\cap \Pi)\right]$$
and 
$$c(\Pi) = \frac{1}{\mathrm{Area}(\Pi)}\left[\mathrm{Area}(\Pi\setminus \Omega_0)\cdot c(\Pi\setminus \Omega_0) + \mathrm{Area}(\Omega_0\cap \Pi)\cdot c(\Omega_0\cap \Pi)\right],$$
whence, after projecting to the $x$-axis and using (\ref{1.37}), we get 
$$c_a(0) =\langle c(\Omega_0),e_1\rangle \le \langle c(\Pi),e_1\rangle =-\frac{B}{2}.
$$

\begin{figure}[ht]
	\centering
	\includegraphics[height=4in, clip=true, viewport=0in 0in 5in 5.5in]{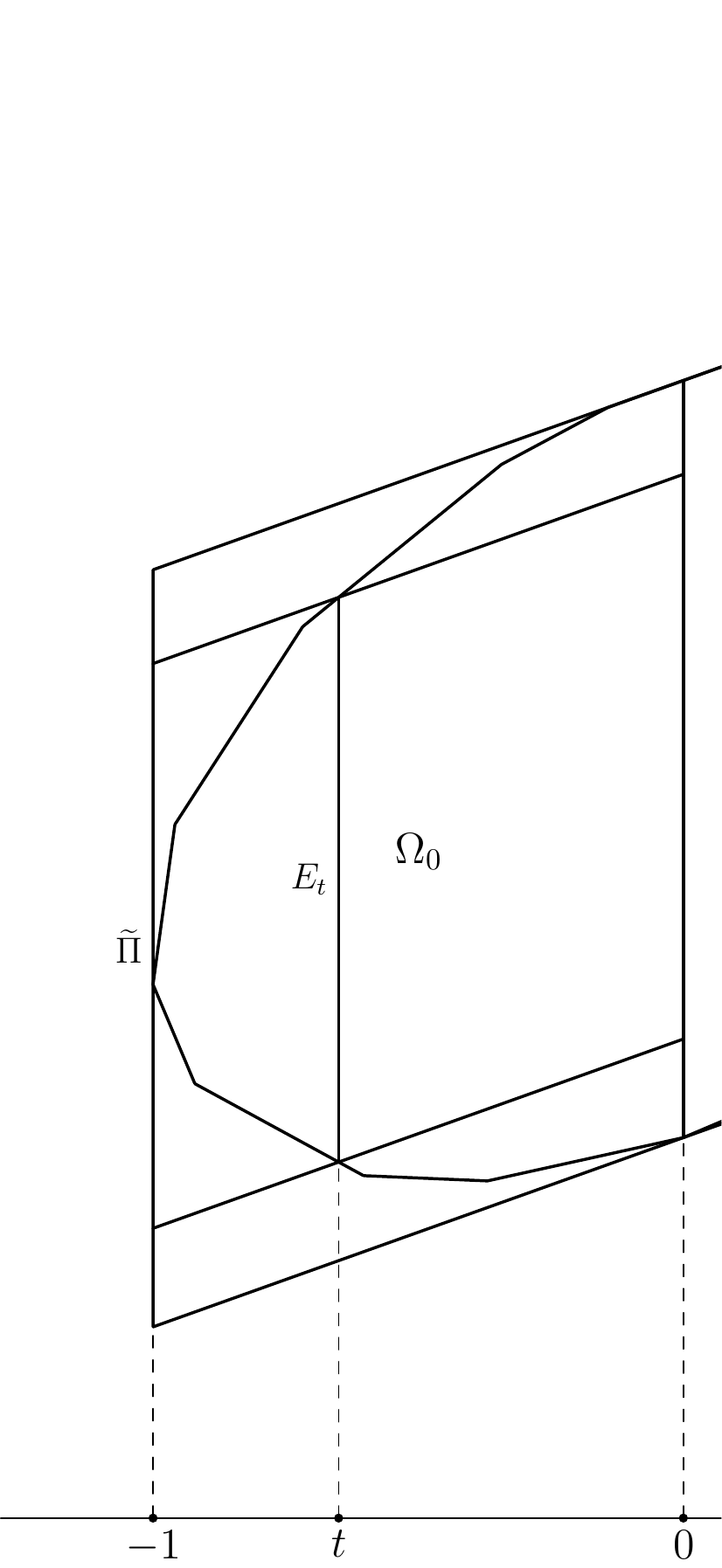}
	\caption{ }
	\label{fig4}
\end{figure}

To get the other inequality, it is enough to compare $\Omega_0$ to the parallelogram $\widetilde\Pi$ bounded by two slanted lines parallel to $L_\pm$ passing through the endpoints of $E_t$ and the vertical lines $x=-1$ and $x=0$ with $t\in (-1,0)$ chosen so that 
$\mathrm{Area}(\widetilde \Pi) = A(0)$  (see Figure \ref{fig4}) 
and to conclude in a similar way that 
$$c_a(0) =\langle c(\Omega_0),e_1\rangle \ge \langle c(\widetilde\Pi),e_1\rangle =-\frac{1}{2}.
$$
\ep

\bl\label{Galya}
Let $\widetilde P(t) = P(t) - \ell(t)$ (the perimeter of $\Omega_t$ {\em without} the right vertical side). 
Then 
$$\frac{A(t)}{A(0)} \le \frac{\widetilde P(t)}{\widetilde P(0)}$$
for all $t\in (-1,0]$.
\el

\bp

\begin{figure}[ht]
 	\centering
 	\includegraphics[height=5in, clip=true, viewport=0in 0in 5in 5.5in]{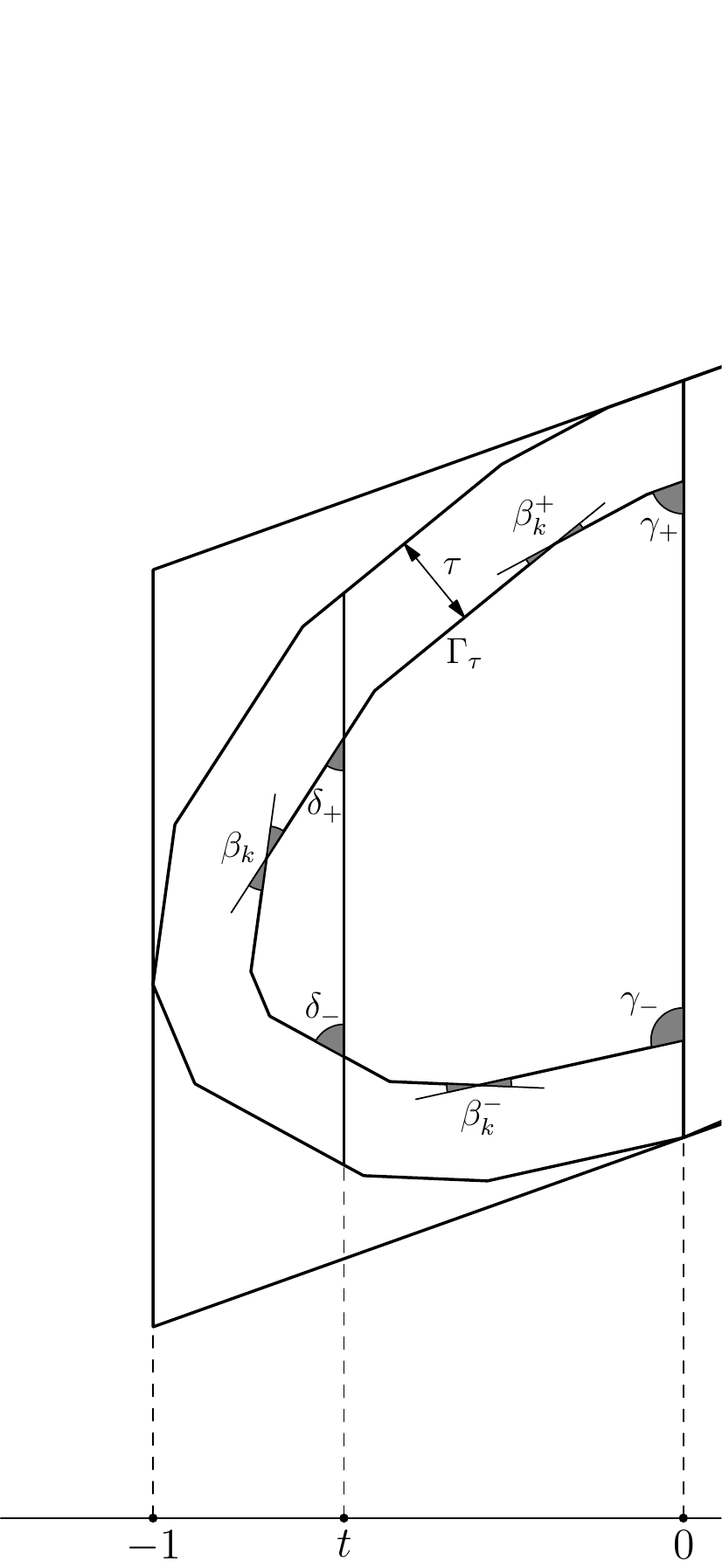}
 	\caption{ }
 	\label{fig5}
 \end{figure}

 Let $\Gamma_\tau$ be the equidistant broken line to $\partial \Omega_0\setminus E_0$ that lies inside $\Omega_0$ and is at distance $\tau$ from $\partial \Omega_0\setminus E_0$. Denote the angles at the vertices of $\Gamma_\tau$ by $\beta_k^-$, $\beta_k$, $\beta_k^+$, where $\beta_k$ lie to the left of the line $x=t$, and $\beta_k^-$,  $\beta_k^+$ between the lines $x=t$ and $x=0$, as shown in Figure \ref{fig5}.

 \noindent
 Also denote by $\gamma_-$ and $\gamma_+$ the   angles between $\Gamma_\tau$ and $x=0$, and by $\delta_-$ and $\delta_+$ the   angles between $\Gamma_\tau$ and $x=t$, as in Figure \ref{fig5}.
Let $\mathfrak{p}_\tau$ be  the part of the perimeter of $\Gamma_\tau$ that lies to the left of the line $x=t$ and 
$\mathfrak{P}_\tau$ the part of the perimeter of $\Gamma_\tau$ that lies to the left of $x=0$. Observe that  the derivatives of $\mathfrak{P}_\tau$ and $\mathfrak{p}_\tau$ with respect to $\tau$ are given by
\begin{multline*} \mathfrak{P}'_\tau = -\sum_k 2 \tan \frac{\beta_k^-}{2} -\sum_k 2 \tan  \frac{\beta_k}{2} -\sum_k 2 \tan  \frac{\beta_k^+}{2} \\ - \tan\left(\frac{\pi}{2}-\gamma_-\right)-\tan\left(\frac{\pi}{2}-\gamma_+\right)
\end{multline*}
and 
$$\mathfrak{p}'_\tau =  -\sum_k 2 \tan \frac{\beta_k  }{2}- \tan\left(\frac{\pi}{2}-\delta_-\right)-\tan\left(\frac{\pi}{2}-\delta_+\right).
$$
All angles  are strictly between $0$ and $\pi$. Since  $\gamma_-+\gamma_+\le \pi$ and $\delta_-+\delta_+\le \pi$, we have $\tan\left(\frac{\pi}{2}-\gamma_-\right)+\tan\left(\frac{\pi}{2}-\gamma_+\right)\ge 0$ and $\tan\left(\frac{\pi}{2}-\delta_-\right)+\tan\left(\frac{\pi}{2}-\delta_+\right)\ge 0$. Thus $\mathfrak{P}'_\tau\le 0$ and $\mathfrak{p}'_\tau\le 0$.

We will now prove the following inequality. If $-\frac{\pi}{2} <\varphi<\varphi+\psi<\frac{\pi}{2}$, then \begin{equation}\label{tan}
\tan(\varphi+\psi)\ge \tan \varphi + 2\tan \frac{\psi}{2}.
\end{equation}
Indeed,
$$\tan(\varphi+\psi) - \tan\varphi = \int_{\varphi}^{\varphi+\psi} \frac{dx}{\cos^2 x}\ge \int_{-\frac{\psi}{2}}^{\frac{\psi}{2}} \frac{dx}{\cos^2 x} = 2 \tan \frac{\psi}{2},$$ where the inequality is a consequence of the two facts that the intervals $[\varphi,\varphi+\psi]$ and $[-\frac{\psi}{2}, \frac{\psi}{2}]$ have the same length, and that $\frac{1}{\cos ^2x}$ is an even function that decreases  on $(-\frac{\pi}{2},0)$ and increases on $(0,\frac{\pi}{2})$.

If  $-\frac{\pi}{2} <\varphi<\varphi+\sum_k\psi_k<\frac{\pi}{2}$, with $\psi_k>0$, then, applying   inequality (\ref{tan}) repeatedly, we obtain 
\begin{equation}\label{tan1}
\tan(\varphi+\sum_k\psi_k)\ge \tan \varphi + \sum_k2\tan \frac{\psi_k}{2}.
\end{equation}

Combining inequality (\ref{tan1}) with the formulas for $\mathfrak{P}_\tau'$, $\mathfrak{p}_\tau'$, and observing that
$$\frac{\pi}{2}-\delta_-  = \frac{\pi}{2}-\gamma_-  +\sum_k \beta_k^-$$ and
$$\frac{\pi}{2}-\delta_+  = \frac{\pi}{2}-\gamma_+  +\sum_k \beta_k^+,$$
we get
$$\mathfrak{p}_\tau'\le \mathfrak{P}_\tau' \le \frac{\mathfrak{p}_0}{\mathfrak{P}_0}\mathfrak{P}_\tau',$$
and, therefore, 
$$\mathfrak{p}_\tau\le   \frac{\mathfrak{p}_0}{\mathfrak{P}_0}\mathfrak{P}_\tau.$$
Let $\mathfrak{s}=\sup\{\tau: \mathfrak{p}_\tau >0 \}$.
Integrating the above inequality with respect to $\tau$, we obtain
$$A(t)=\int_0^{\mathfrak{s}} \mathfrak{p}_\tau d\tau\le \frac{\mathfrak{p}_0}{\mathfrak{P}_0} \int_0^{\mathfrak{s}} \mathfrak{P}_\tau d\tau \le \frac{\mathfrak{p}_0}{\mathfrak{P}_0} A(0).$$
Since $\mathfrak{p}_0=\widetilde P(t)$ and  $\mathfrak{P}_0=\widetilde P(0)$, the desired inequality follows.
\ep

\bl\label{LemmaTrapezoid} Let $t\in (-1,0)$ and $\ell'=\ell(t)$. Then
$$\frac{A(t)}{A(0)} \le (1+t) \frac{(1+t) \ell -(1-t)\ell'}{(1+2t)\ell - \ell'}.$$

\el

\begin{figure}[ht]
	\centering
	\includegraphics[height=4.5in, clip=true, viewport=0in 0in 5in 6.8in]{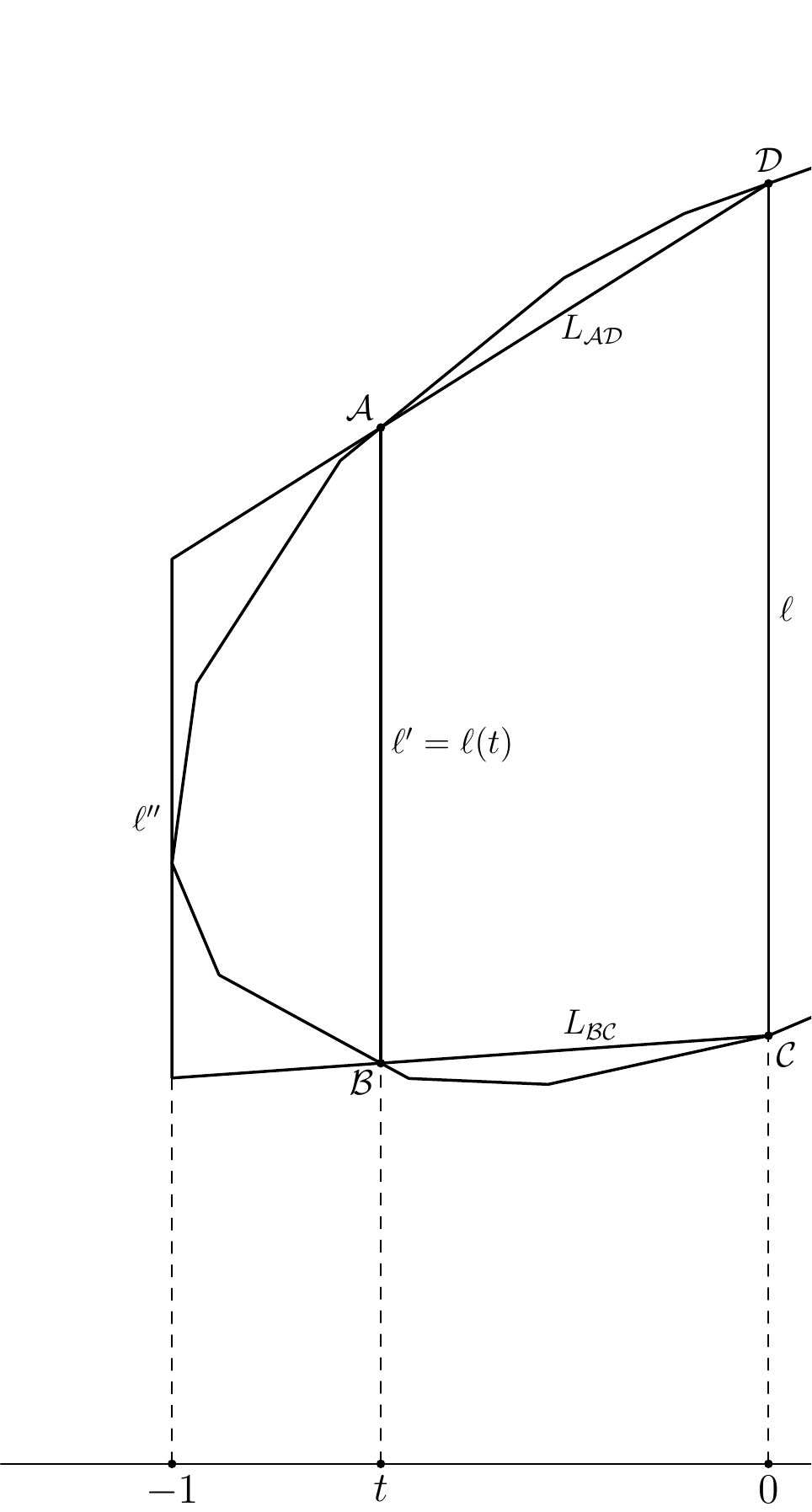}
	\caption{ }
	\label{figlemma3}
\end{figure}

\bp
Let $\mathcal A$, $\mathcal B$, $\mathcal C$, and $\mathcal D$ be the points of intersection of the boundary of $\Omega$ and the lines $x=t$ and $x=0$; see Figure \ref{figlemma3}.

\noindent
Denote by $L_{\mathcal A \mathcal D}$  the line  passing through $\mathcal A$ and $\mathcal D$, and by $L_{\mathcal B \mathcal C}$  the line  passing through $\mathcal B$ and $\mathcal C$. 
Consider the trapezoid $\mathcal T$ bounded by the lines $x=-1$, $x=0$, $L_{\mathcal A \mathcal D}$, and $L_{\mathcal B \mathcal C}$. Denote by $\ell''$ the side of this trapezoid that lies on the line $x=-1$. Let $S(t)$  be the area of $\mathcal T\cap \{x\le t\}$, and $S(0)$  the total area of $\mathcal T$. Then, by inclusion, we have $A(t)\le S(t)$ and $A(0)-A(t)\ge S(0)-S(t)$. Hence 
$$\frac{A(0)-A(t)}{A(t)}\ge \frac{S(0)-S(t)}{S(t)}.$$
This gives
$$\frac{A(t)}{A(0)}\le \frac{S(t)}{S(0)} = (1+t)\frac{\ell'+\ell''}{\ell+\ell''}.$$
Observe that $\ell'=\ell''+(1+t)(\ell - \ell'')$, i.e., $\ell'' = \frac{(1+t)\ell-\ell'}{t}$. Substituting this expression for  $\ell''$ into the above inequality, we obtain the statement of the lemma.
\ep

\begin{figure}[ht]
	\centering
	\includegraphics[height=4in]{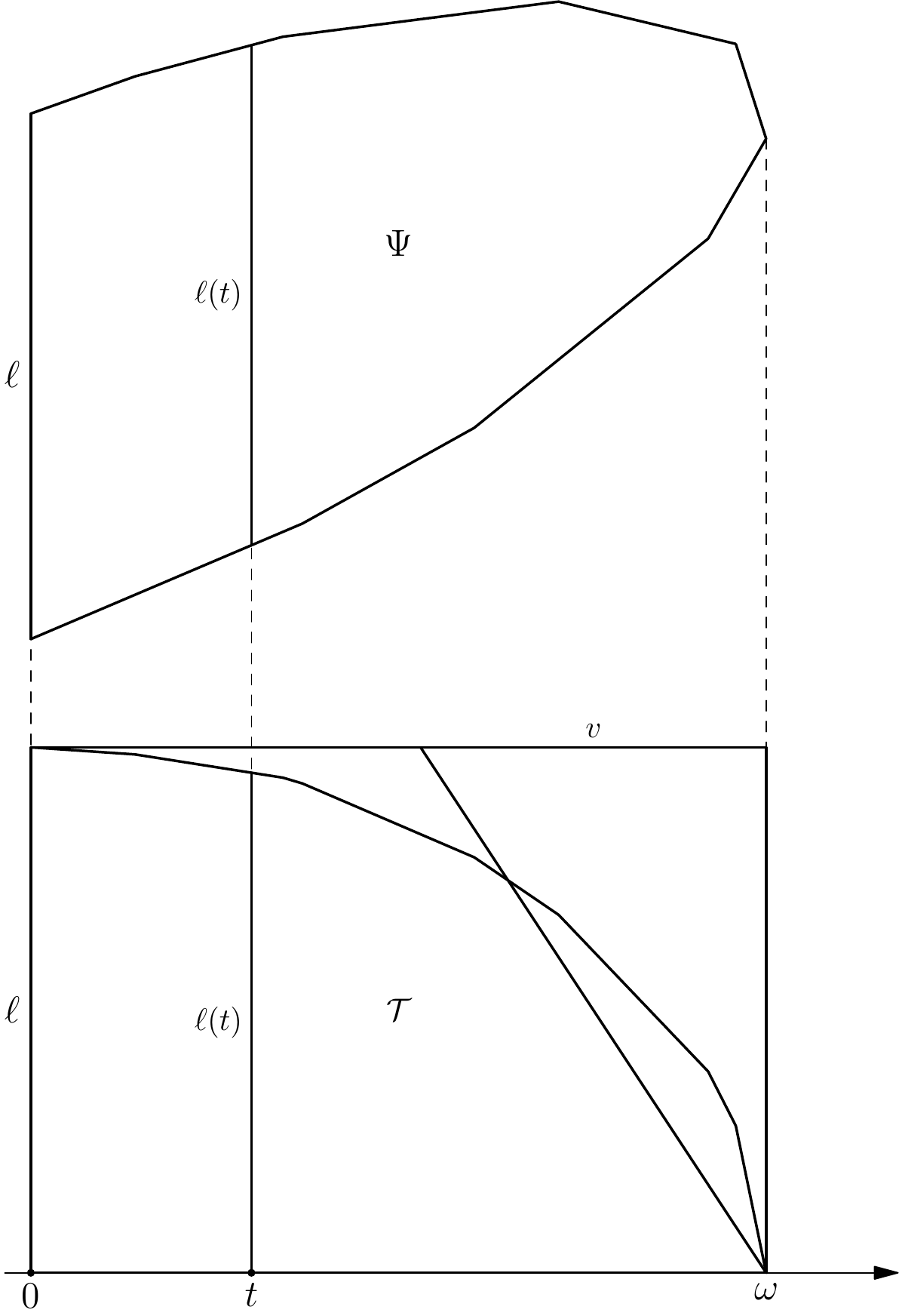}
	\caption{ }
	\label{figlemma4}
\end{figure}

\bl\label{fedsm}
Let $ {c}_\Psi$ be the $x$-coordinate of the center of mass of the body $\Psi=\Omega\setminus \Omega_0$ shown on Figure \ref{figlemma4}. Put $u=2 \frac{\mathrm{Area}(\Psi)}{\ell \omega}\in [1,2]$. We have 
$$
{c}_\Psi\ge\omega\frac{u^2-u+1}{3u}.
$$
\el

\bp
The $x$-coordinate of the center of mass of $\Psi$ does not change if we move the intersection $E_t$ of $\Psi$ with every vertical line $x=t$, $t\in [0,\omega]$, in the vertical direction so that its lower end lies on the $x$-axis. The resulting domain will stay convex and be contained in the rectangle $[0,\omega]\times [0,\ell]$. Construct the slanted line passing through the lower right corner of this rectangle so that the area of the resulting
trapezoid ${\mathcal T}$ is  the same  as of $\Psi$ (see Figure~\ref{figlemma4}).  

\begin{figure}[ht]
	\centering
	\includegraphics[height=4in, clip=true, viewport=0in 0in 6in 6in]{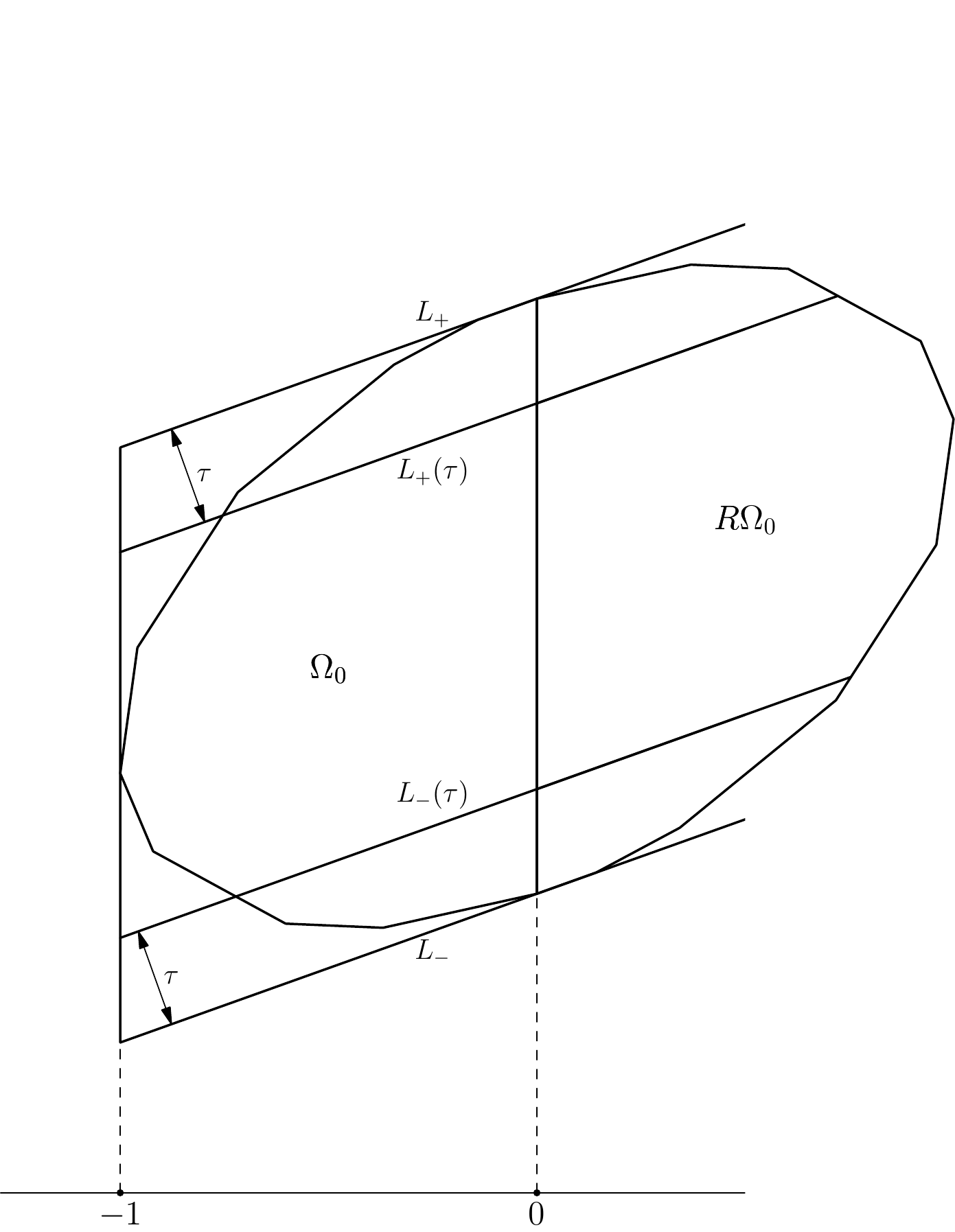}
	\caption{ }
	\label{figlemma5}
\end{figure}

Observe that the $x$-coordinate of  $c(\mathcal T)$ is
at most ${c}_\Psi$. Denote by $v$ the length of the upper side of the triangle  $([0,\omega]\times [0,\ell])\setminus \mathcal T$. Since the area of the trapezoid is $\frac{(\omega-v)+\omega}{2}\ell$, we have  $v=(2-u)\omega$. Evaluating the $x$-coordinate of the center of mass of the trapezoid as a linear combination of the centers of mass of the rectangle and the triangle, we see that
$$
{c}_\Psi\ge\frac{\frac{\omega}{2}\omega\ell-(\omega-\frac{1}{3}v)\cdot\frac{1}{2}v\ell}{\frac{u\ell\omega}{2}}=\omega\frac{u^2-u+1}{3u}.
$$
\ep
\bl\label{dfedsm1}
We have
$$
A(0)\le \frac{P(0)-\ell}{2\alpha}\cdot\ell.
$$
\el
\bp
Let $L_-(\tau)$ and $L_+(\tau)$ be the two  lines intersecting the interior of $\Omega$ that are parallel to $L_-$ and $L_+$ and have distance $\tau$ from $L_-$ and $L_+$ respectively (see Figure~\ref{figlemma5}). 
Let $\widetilde{\Omega}=\Omega_0\cup R\Omega_0$ where  $R\Omega_0$ is the refection of $\Omega_0$  about the midpoint of $E_0$. Observe that
$$
\textrm{length}(L_+(\tau)\cap\Omega_0))+\textrm{length}(L_-(\tau)\cap\Omega_0))=
\textrm{length}(L_+(\tau)\cap\widetilde{\Omega})$$
$$
\le\,\frac{1}{2}\,\textrm{length}(\partial\widetilde{\Omega})=P(0)-\ell.
$$
Then, since the width of the strip bounded by   $L_+$ and $L_-$ is   $2h=\frac{\ell}{\alpha}$, we conclude that 
\begin{gather*}
A(0)=\int\limits_0^h(\textrm{length}(L_+(\tau)\cap\Omega_0)+\textrm{length}(L_-(\tau)\cap\Omega_0))\,dt
\\
\le \int\limits_0^h(P(0)-\ell)\,dt=(P(0)-\ell)h=\frac{P(0)-\ell}{2\alpha}\cdot\ell.
\end{gather*}
\ep

\begin{figure}[ht]
 	\centering
 	\includegraphics[height=5in]{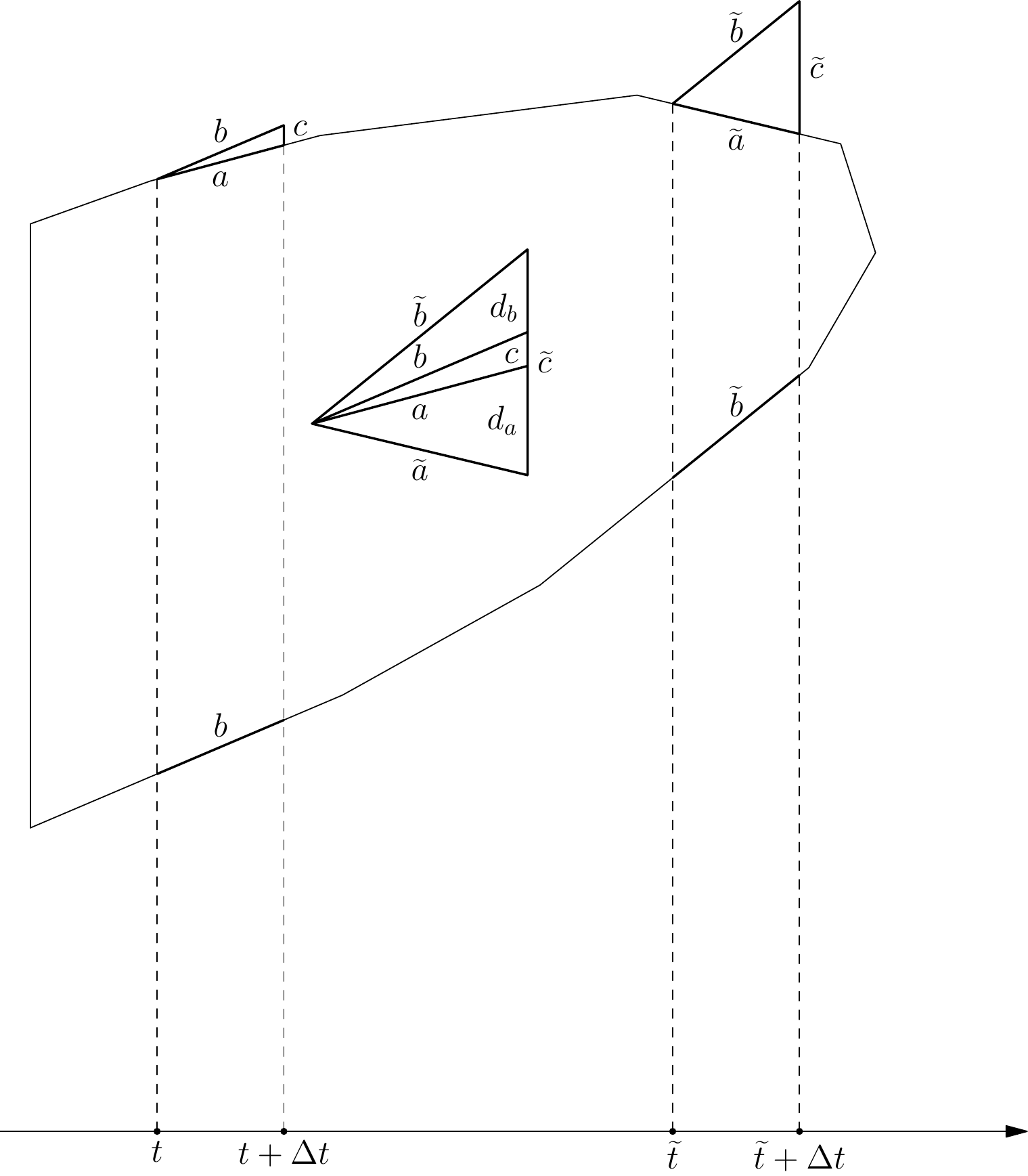}
 	\caption{ }
 	\label{figlemma6}
 \end{figure}

\bl\label{P(t)concave}
$P'(t)$ is a decreasing  function of $t$.
\el
\bp
 We will prove the lemma for $t>0$, which is the only case that will be needed later. However, an easy modification of the argument below will also yield the case $t< 0$.   Note that $P'(t)$ exists for all $t$ except finitely many values. Let $t$, $\widetilde t$, $\Delta t$ be such that $t<t+\Delta t<\widetilde t <\widetilde t +\Delta t$ and there is no vertex with the $x$-coordinate in $(t,t+\Delta t)\cup (\widetilde t, \widetilde t +\Delta t)$.   Let $a$, $b$, $\widetilde a$, $\widetilde b$ be the lengths shown on Figure \ref{figlemma6},  $c=\ell(t)- \ell(t+\Delta t)$, and $\widetilde c= \ell(\widetilde t)- \ell(\widetilde t+\Delta t)$. 
 
 Then $P'(t)\Delta t=a+b-c$ and $P'(\,{\widetilde t}\,)\Delta t=\widetilde a+\widetilde b-\widetilde c$. We need to show that $$a+b-c\ge \widetilde a+\widetilde b-\widetilde c.$$
 This  follows from the triangle inequality applied twice:   
 $\widetilde a\le a+d_a$ and $\widetilde b\le b +d_b$ (see  Figure~\ref{figlemma6}), and the equality $d_a+d_b=\widetilde c - c$. 
 \ep

\section{Proof of inequality  (\ref{**}).}\label{Sec5} Denote $$p(t) = \frac{\widetilde P(t)}{  P(0)}\quad \mbox{and} \quad a(t)= \frac{A(t)}{A(0)}.$$ We have 

$$
c_p(0)-c_a(0) = \int_{-1}^0 t \, d\left( p(t) - a(t) \right) =\int_{-1}^0  \left(  a(t) - p(t)\right)\,dt.$$
We will now estimate $ a(t) - p(t)$ in three different ways. 
First, by Lemma \ref{Galya},

\begin{equation}\label{3ineq1} a(t) - p(t)   \le \frac{\widetilde P(t)}{\widetilde P(0)} - \frac{\widetilde P(t)}{\widetilde P(0) +\ell} = \frac{\widetilde P(t)\ell}{(\widetilde P(0)+\ell)\widetilde P(0)}.
\end{equation}
Secondly, by Lemma \ref{LemmaTrapezoid},
\begin{equation}\label{a(t)} a(t) \le (1+t) \frac{2\ell'-(1+t) (\ell' +\ell)}{\ell+\ell'-2(1+t)\ell}. 
\end{equation}
Since the mapping $z\mapsto  \frac{2z-(1+t) (z +\ell)}{\ell+z-2(1+t)\ell}$ is fractional-linear in $z$, it is monotone on any ray   $(z_0, +\infty)$  on which the denominator stays positive. Observing that its value 1 at $z=\ell$ is less than its value $1-t$ at $z=+\infty$, we conclude that it is increasing in $z$. 
The last observation together with the inequality $\ell'\le \widetilde P(t)$ implies
\begin{equation}\label{3ineq2}a(t)-p(t)\le (1+t) \frac{2\widetilde P(t)-(1+t) (\widetilde P(t) +\ell)}{\ell+\widetilde P(t)-2(1+t)\ell} - \frac{\widetilde P(t)}{\widetilde P(0)+\ell}.
\end{equation}
Finally, using  (\ref{a(t)}) with  $\ell'$ replaced by $ \ell$ , 
 we get 
 $ a(t)\le   1+t,$ i.e., \begin{equation}\label{3ineq3} a(t)-p(t)\le 1+t - \frac{\widetilde P(t)}{\widetilde P(0)+\ell}.
 \end{equation}
Putting 
$$\lambda =  \frac{\ell}{\widetilde P(0)}, \quad \mu = \frac{\widetilde P(t)}{\widetilde P(0)}\in [0,1],$$  we can rewrite inequalities (\ref{3ineq1}), (\ref{3ineq2}), and (\ref{3ineq3})  as 
\begin{equation}\label{fed1}
a(t)\!-p(t)\!\le \!
\min\!\left[ (1+t)\!\cdot\!\min\!\!\left(1,\frac{ 2\mu - (1+t) (\mu+\lambda)}{\lambda+\mu-2(1+t)\lambda}\!\right)\!-\frac{\mu}{1+\lambda}, \frac{\lambda\mu}{1+\lambda}\!\right]\!.
\end{equation}
We will use the above inequality to estimate the integral
$$I=\int_{-1}^0  \left(  a(t) - p(t)\right) \, dt=\int_{-1}^{\lambda-1}  \left(  a(t) - p(t)\right) \, dt+\int_{\lambda-1}^0 \left(  a(t) - p(t)\right) \, dt.$$
For $t\in [-1,\lambda -1]$, inequality (\ref{fed1}) yields
\begin{align*}
a(t) - p(t)  &\le (1+t)\cdot \frac{ 2\mu - (1+t) (\mu+\lambda)}{\lambda+\mu-2(1+t)\lambda
} -\frac{\mu}{1+\lambda}\\
&= (1+t)\cdot \frac{ \mu(1-t) + \lambda (1-t) -2 \lambda (1+t) (1-t) -2 \lambda t^2  }{\lambda+\mu-2(1+t)\lambda
} -\frac{\mu}{1+\lambda}\\
&=  (1+t)\cdot\left(1-t - \frac{ 2\lambda t^2}{\lambda+\mu-2(1+t)\lambda
}\right) -\frac{\mu}{1+\lambda}\\
&= 1-t^2  - \frac{ 2\lambda t^2(1+t)}{\lambda+\mu-2(1+t)\lambda
}  -\frac{\lambda+\mu-2(1+t)\lambda }{1+\lambda}-\frac{\lambda (1+2t)}{1+\lambda}.
\end{align*}
Using the AM-GM  inequality to estimate 
two middle terms, we  obtain 
$$
a(t) - p(t)  \le  1-t^2  + 2t\sqrt{\frac{ 2\lambda (1+t)}{1+\lambda
}}  -\frac{\lambda (1+2t)}{1+\lambda}.$$
On the other hand, for $t\in [\lambda -1,0]$, we will use (\ref{fed1}) in the form
$$a(t) - p(t)  \le \min\left[ 1+t -\frac{\mu}{1+\lambda}, \frac{\lambda\mu}{1+\lambda}\right].$$
Note that when $t\le \mu -1$, the latter minimum is given by the first function, and when $t\ge \mu -1$, it is given by the second function. In both cases we have 
$$a(t) - p(t)  \le   \frac{\lambda (1+t)}{1+\lambda} .$$
Now we have
\begin{align*}&I=\int_{-1}^{\lambda-1}  \left(  a(t) - p(t)\right) \, dt+\int_{\lambda-1}^0 \left(  a(t) - p(t)\right) \, dt\\
&\le \int_{-1}^{\lambda-1}  \left(   1-t^2  + 2t\sqrt{\frac{ 2\lambda (1+t)}{1+\lambda
}}  -\frac{\lambda (1+2t)}{1+\lambda}\right) \, dt+\int_{\lambda-1}^0   \frac{\lambda (1+t)}{1+\lambda} \, dt.
\end{align*}
Making the change of variable $t=\tau-1$ in the first integral, we get 
\begin{align*}&I\le \int_{0}^{\lambda}  \left(   2\tau-\tau^2  + 2\tau \sqrt{\frac{ 2\lambda \tau}{1+\lambda
}} - 2 \sqrt{\frac{ 2\lambda \tau}{1+\lambda
}} -\frac{\lambda (2\tau -1 )}{1+\lambda}\right) \, dt+  \frac{\lambda (1+t)^2}{2(1+\lambda)}  \Big|_{t=\lambda-1}^{t=0} \\
&= \int_{0}^{\lambda}  \left(  \frac{\lambda }{1+\lambda}+  \frac{2}{1+\lambda}\tau -\tau^2  + 2\sqrt{\frac{ 2\lambda  }{1+\lambda
}} \tau^{3/2} - 2 \sqrt{\frac{ 2\lambda  }{1+\lambda
}} \tau^{1/2} \right) \, dt+\frac12\lambda(1-\lambda)\\
&=   \frac{\lambda^2 }{1+\lambda}+  \frac{\lambda^2}{1+\lambda}-\frac{\lambda^3}{3}  + \frac45\sqrt{\frac{ 2\lambda  }{1+\lambda
}} \lambda^{5/2} - \frac43 \sqrt{\frac{ 2\lambda  }{1+\lambda
}} \lambda^{3/2} +\frac12\lambda(1-\lambda) \\
&= \frac{2\lambda^2 }{1+\lambda}-\frac13 \lambda^3+4 \sqrt{\frac{ 2   }{1+\lambda }}\left( \frac15 \lambda^{3} -  \frac13 \lambda^{2}\right)+\frac12\lambda(1-\lambda).
\end{align*}
Note that $\sqrt{\frac{ 2  }{1+\lambda }}$ is a convex function of $\lambda$, so its graph lies above its tangent line at $\lambda=1$, i.e., 
$\sqrt{\frac{ 2  }{1+\lambda }} \ge 1-\frac{1}{4}(\lambda-1)=\frac14(5-\lambda)$.
Since $\frac15 \lambda^{3} -  \frac13 \lambda^{2}<0$ on $[0,1]$, we conclude that
$$I\le -\frac13 \lambda^3+(5-\lambda)\left( \frac15 \lambda^{3} -  \frac13 \lambda^{2}\right)+\frac{2\lambda^2}{1+\lambda}+\frac12\lambda(1-\lambda).$$
To show that the latter does not exceed $\frac16$, we will multiply it by $30(1+\lambda)$ and compare  to $ 5 (1+\lambda)$. We get
$$-10(1+\lambda) \lambda^3+(1+\lambda)(5-\lambda)\left( 6 \lambda^{3} -  10 \lambda^{2}\right)+60 \lambda^2 +15\lambda(1-\lambda^2)\le 5 (1+\lambda),$$
which is equivalent to
\begin{multline*}-10  \lambda^3  -10  \lambda^4 +(5+4\lambda -\lambda^2)\left( 6 \lambda^{3} -  10 \lambda^{2}\right)\\ +60 \lambda^2 +15\lambda-15\lambda^3 
\le 5+5\lambda.
\end{multline*}
Simplifying it further, we get
\begin{multline*}-10  \lambda^3  -10  \lambda^4+30\lambda^3-50\lambda^2+24\lambda^4-40\lambda^3-6\lambda^5 +10\lambda^4 \\ +60 \lambda^2 +15\lambda-15\lambda^3
\le 5+5\lambda.
\end{multline*}
Finally, the inequality that we need to prove becomes 
$$5-10\lambda -10\lambda^2+35\lambda^3-24\lambda^4+6\lambda^5\ge 0, \quad \mbox{for } \lambda\in[0,1].$$
This inequality follows from the fact that the above polynomial can be expressed in the form
\begin{align*}&5-10\lambda -10\lambda^2+35\lambda^3-24\lambda^4+6\lambda^5\\
&\qquad=(5-16\lambda +13\lambda^2)+(\lambda^2-\lambda^3)+(6\lambda-24\lambda^2+36\lambda^3-24\lambda^4+6\lambda^5)\\
&\qquad=(5-16\lambda +13\lambda^2)+(\lambda^2-\lambda^3)+6\lambda(1-\lambda)^4.
\end{align*}

\section{Proof of inequality (\ref{***})}

We start with an upper bound for $c_p(\omega)$. 
Observe that
\begin{align*}P(t)c_p(t)& = \int_{-1}^t s\, d(P(s)-\ell(s)) + t \ell(t)\\
&=t(P(t)-\ell(t))+\int_{-1}^t (\ell(s)-P(s))ds +t\ell(t)\\
&=t P(t) - \int_{-1}^t P(s)ds +A(t).
\end{align*}
Using the above formula with $t=0$ and $t=\omega$, we get
\begin{align*} 
P(0)c_p(0)&=-\int\limits_{-1}^0P(s)ds+A(0),
\\
P(\omega)c_p(\omega)&=P(\omega)\omega-\int\limits_{-1}^\omega P(s)ds+A(\omega).
\end{align*}
This gives
\begin{align*} 
P(\omega)c_p(\omega)&=P(0)c_p(0)+P(\omega)\omega-\int\limits_{0}^\omega P(s)ds
+A(\omega)-A(0)\\
&\le P(0)c_p(0)+P(\omega)\omega-\frac{P(0)+P(\omega)}{2} \omega+A(\omega)-A(0),
\end{align*}
where in the last estimate we used the concavity of $P(t)$; see Lemma \ref{P(t)concave}.
Dividing both parts by $P(\omega)$ we obtain
\begin{align*}
c_p(\omega)&\le \frac{P(0)}{P(\omega)}c_p(0)+\omega-\frac{P(0)+P(\omega)}{2P(\omega)} \omega+\frac{A(\omega)-A(0)}{P(\omega)}\\
&=
\frac{\omega}{2}+\frac{P(0)}{P(\omega)}c_p(0)-\frac{P(0)\omega-u\ell\omega}{2P(\omega)},
\end{align*}
where we chose  $u\in [1,2]$ so that $A(\omega)-A(0)=\frac{u\ell\omega}{2}$ as in Lemma \ref{fedsm}.
Observe that $$
\mathrm{length}(\partial \Omega\cap\{x\ge 0\})\le \ell+2\omega \alpha,
$$ 
and therefore 
$$P(\omega)=P(0)-\ell+\mathrm{length}(\partial \Omega \cap\{x\ge 0\})\le P(0)+2\omega \alpha.$$
Since $c_p(0)<0$, and since the inequality $P(0)\ge 2\ell$ yields  $
P(0)\omega\ge u\ell\omega$, we obtain
$$
c_p(\omega)\le\frac{\omega}{2}+\frac{P(0)}{P(0)+2\omega \alpha}c_p(0)-\frac{P(0)\omega-u\ell\omega}{2(P(0)+2\omega \alpha)}.
$$
Setting  $c=-c_a(0)$ and using inequality (\ref{**}), we see that the above inequality  yields
\begin{equation}\label{fedsm1}
c_p(\omega)\le\frac{\omega}{2}+\frac{P(0)}{P(0)+2\omega \alpha}\Big(\frac{1}{6}-c\Big)-\frac{P(0)\omega-u\ell\omega}{2(P(0)+2\omega \alpha)}.
\end{equation}
Now we evaluate $c_a(\omega)$ as 
$$
c_a(\omega)=\frac{A(0)c_a(0)+\frac{u\ell\omega}{2}{c}_\Psi}{A(0)+\frac{u\ell\omega}{2}},
$$
where ${c}_\Psi$ is the $x$-coordinate of the center of mass of $\Psi=\Omega\cap \{x\ge 0\}$.
Using Lemma~\ref{fedsm},  we obtain
 \begin{equation}\label{fedsm2}
c_a(\omega)\ge\frac{-A(0)c+\frac{1}{6}\ell\omega^2(u^2-u+1)}{A(0)+\frac{u\ell\omega}{2}} .
 \end{equation}

Combining (\ref{fedsm1}) and (\ref{fedsm2}), we see that (\ref{***}) would follow from 
\begin{multline}\label{eqnn}
\frac{\omega}{2}+\frac{P(0)}{P(0)+2\omega \alpha}\Big(\frac{1}{6}-c\Big)-\frac{P(0)\omega-u\ell\omega}{2(P(0)+2\omega \alpha)}
\\
-\,\frac{-A(0)c+\frac{1}{6}\ell\omega^2(u^2-u+1)}{A(0)+\frac{u\ell\omega}{2}}
\le  \frac{1+\omega}{6}.
\end{multline}
We now move all the terms to one side and denote  $P=P(0)$, $B=\frac{A(0)}{\ell} $. Then (\ref{eqnn}) rewrites as
$$
\frac{1}{6}-\frac{\omega}{3}+\frac{P(c-\frac{1}{6})}{P+2\omega \alpha}+\frac{P\omega-u\ell\omega}{2(P+2\omega \alpha)}+
\frac{ -Bc+\frac{1}{6}\omega^2(u^2-u+1)}{B+\frac{u\omega}{2}}\ge 0,
$$
or, equivalently, as
$$
\frac{1}{6}-\frac{\omega}{3} +c-\frac{1}{6}+\omega\frac{-2\alpha (c-\frac{1}{6}) +\frac{P-u\ell}{2}  }{P+2\omega\alpha}            -c+\omega\frac{\frac{uc}{2}+\frac{1}{6}\omega(u^2-u+1)}{B+\frac{u\omega}{2}} \ge 0.
$$
Canceling the constant terms and dividing by $\omega$, we see that it suffices  to prove that
$$
\frac{\frac{uc}{2}+\frac{1}{6}\omega(u^2-u+1)}{B+\frac{u\omega}{2}}+\frac{-2\alpha (c-\frac{1}{6}) +\frac{P-u\ell}{2}  }{P+2\omega\alpha}   -\frac{1}{3}\ge 0.
$$
Now  set $s=\frac{2\alpha}{P}, \lambda = \frac{2\ell}{P}$. It is clear that $0<\lambda <1$ and $s>0$. To see that $s<1$, draw the line
parallel to $L_\pm$ through the single point in $E_{-1}$. This line
divides $\partial \Omega_0$ into two parts of length greater than
$\alpha$.
Thus, 
$s,\lambda \in (0,1)$ and    we can rewrite the previous inequality as
$$
\frac{\frac{uc}{2}+\frac{1}{6}\omega(u^2-u+1)}{B+\frac{u\omega}{2}}+\frac{-sP(c-\frac{1}{6}) +\frac{P}{2}-\frac{u\lambda P}{4}  }{P+\omega sP}   -\frac{1}{3}\ge 0.
$$
Next, we set $s\omega=\rho$ and $sB=b$. Canceling $P$ in the second summand, we obtain
$$
\frac{\frac{suc}{2}+\frac{1}{6}\rho(u^2-u+1)}{b+\frac{u\rho}{2}}+\frac{-s(c-\frac{1}{6}) +\frac{1}{2}-\frac{u\lambda }{4}  }{1+\rho}   -\frac{1}{3}\ge 0.
$$
This is equivalent to
$$
\frac{1}{3u}(u^2-u+1)+\frac{\frac{suc}{2}-\frac{b}{3u}(u^2-u+1)}{b+\frac{u\rho}{2}}+\frac{-s(c-\frac{1}{6}) +\frac{1}{2}-\frac{u\lambda }{4}  }{1+\rho}   -\frac{1}{3}\ge 0.
$$
It is enough to show the stronger inequality
$$
\frac{u-2+\frac{1}{u}}{3(1+\rho)}+\frac{\frac{suc}{2}-\frac{b}{3u}(u^2-u+1)}{b+\frac{u\rho}{2}}+\frac{-s(c-\frac{1}{6}) +\frac{1}{2}-\frac{u\lambda }{4}  }{1+\rho}  \ge 0.
$$
After  multiplying by the denominators, the left hand side  becomes a linear function in $\rho$. To show that this  function is non-negative on $[0,+\infty)$, we need to prove that its coefficients are non-negative. One of them is evaluated by setting $\rho=0$ and the other one (up to a factor $\frac{u}{2}$) is obtained by multiplying  the previous inequality by $\rho$ and sending $\rho$ to infinity.
In other words, we have to check the following two inequalities:
\begin{equation}\label{fedsm3}
\frac{u-2+\frac{1}{u}}{3}+\frac{\frac{suc}{2}-\frac{b}{3u}(u^2-u+1)}{b}-s\Big(c-\frac{1}{6}\Big) +\frac{1}{2}-\frac{u\lambda }{4}   \ge 0,
\end{equation}
and
\begin{equation}\label{fedsm4}
\frac{u-2+\frac{1}{u}}{3}+sc-\frac{2b(u^2-u+1)}{3u^2}-s\Big(c-\frac{1}{6}\Big) +\frac{1}{2}-\frac{u\lambda }{4} \ge 0.
\end{equation}
From now on, we will forget about the geometric meanings of $c$, $u$, $s$, $\lambda$, $B$, $b$ and will treat them as free variables about which we only assume the restrictions $u\in [1,2]$, $s,\lambda\in (0,1]$, $B\in [\frac{1}{2},\min(1,\frac{1-\frac{\lambda}{2}}{s})]$, $c\in [\frac{B}{2},\frac{1}{2}]$, $b=sB$.

 First, we check (\ref{fedsm3}). Substituting $b=sB$, we can rewrite it  as
\begin{equation}\label{fedsm5}
\frac{1}{6}+\frac{uc}{2B}-s\Big(c-\frac{1}{6}\Big)-\frac{u\lambda}{4}\ge 0.
\end{equation}
By Lemma \ref{dfedsm1} 
	we have $B\le \frac{1-\frac{\lambda}{2}}{s}$, and so $\lambda\le 2(1-sB)$.  Also,  $\lambda\le 1$. Hence, 
	$$
	\frac{u\lambda}{4}\le \frac{u}{2}\min\Big(\frac{1}{2}, 1-sB\Big)=\frac{u}{2}\Big(1-\max\Big(sB, \frac{1}{2}\Big)\Big)
	$$
	and  (\ref{fedsm5})  will follow from
	\begin{equation}\label{fedsm333}
\frac{1}{6}+\frac{uc}{2B}-s\Big(c-\frac{1}{6}\Big)-\frac{u}{2}+\frac{u}{2}\max\Big(sB, \frac{1}{2}\Big)\ge 0.
\end{equation}

 By Lemma \ref{dfedsm22}, $\frac{B}{2}\le c\le \frac{1}{2}$.
Since the left hand side of (\ref{fedsm333})  is linear in $c$, it is enough to check it for $c=\frac{1}{2}$ and 
$c=\frac{B}{2}$. 

Let $c=\frac{1}{2}$. It is enough to show that
\begin{equation}\label{fedsm111}
\frac{1}{6}+\frac{u}{4B}-\frac{s}{3}-\frac{u}{2}+\frac{u}{2}\max\Big(sB, \frac{1}{2}\Big)\ge 0.
\end{equation}
If $s\le \frac{1}{2B}$, the left-hand side of (\ref{fedsm111}) is a decreasing linear function of $s$, and so its minimum is achieved at $s=\frac{1}{2B}$. Thus, it is enough to check (\ref{fedsm111}) only when $s\ge \frac{1}{2B}$, i.e., that
\begin{equation}\label{fedsm1111}
\frac{1}{6}+\frac{u}{4B}-\frac{s}{3}-\frac{u}{2}+\frac{u}{2}sB  \ge 0.
\end{equation}
The left-hand side  of (\ref{fedsm1111}) is linear in $s$,  and so we need to check its sign at the endpoints of the interval $ [\frac{1}{2B},1]$.  If $s=1$, we apply the AM-GM inequality and then use $u\ge 1$  to get
$$
\frac{1}{6}+\frac{u}{4B}-\frac{1}{3}-\frac{u}{2}+\frac{uB}{2}\ge -\frac{1}{6}-\frac{u}{2}+
\frac{u}{\sqrt{2}}\ge \frac{u}{6}-\frac{1}{6}\ge 0.
$$
If
$s =\frac{1}{2B}$, then (\ref{fedsm1111}) reads as
$$
\frac{1}{6}+\frac{u}{4B}-\frac{1}{6B}-\frac{u}{2}+\frac{u}{4}=\Big(\frac{1}{B}-1\Big)\Big(\frac{u}{4}-\frac{1}{6}\Big)\ge 0,
$$
where the last inequality is due to the fact that $B\in[\frac{1}{2},1]$ and $u\ge 1$.

Now let $c=\frac{B}{2}$. Then (\ref{fedsm333}) rewrites as
$$
\frac{1}{6}+\frac{u}{4}-s\Big(\frac{B}{2}-\frac{1}{6}\Big)-\frac{u}{2}+\frac{u}{2}\max\Big(sB, \frac{1}{2}\Big) \ge 0.
$$
Since $\frac{B}{2}-\frac16>0$ and $\frac{Bu}{2}-\left(\frac{B}{2}-\frac16\right)>0$, the  minimum of the left-hand side as a function of $s$ is achieved at $s=\frac{1}{2B}$ and equals
$\frac{1}{12B}-\frac{1}{12}$ which is non-negative because $B\le 1$.

It remains to check (\ref{fedsm4}), i.e., 
\begin{equation}\label{EEq}
\frac{u-2+\frac{1}{u}}{3}-\frac{2b(u^2-u+1)}{3u^2}+\frac{s}{6} +\frac{1}{2}-\frac{u\lambda }{4} \ge 0.
\end{equation}
Recall that $b\le \min(1-\frac{\lambda}{2}, s)$. We will consider two cases, {\it Case 1}: $s\le \frac{1}{2}$,  and {\it Case 2}: $s\ge \frac12$.

{\it Case 1}. Since $s\le \frac{1}{2}$, we have $\min(1-\frac{\lambda}{2}, s)=s$. The left hand side of (\ref{EEq}) decreases when  we replace $\lambda$ and $b$ by their largest values, i.e., $\lambda=1$ and $b=s$, which yields
\begin{align*}
&\frac{u-2+\frac{1}{u}}{3}-\frac{2s(u^2-u+1)}{3u^2}+\frac{s}{6} +\frac{1}{2}-\frac{u }{4} \\
& =
-\frac{1}{6}+\frac{u}{12}+\frac{1}{3u}-s\Big(\frac{1}{2}-\frac{2}{3u}+\frac{2}{3u^2}   \Big)\\
& =-\frac{1}{6}+\frac{u}{12}+\frac{1}{3u}-s\frac{1}{6u^2}\left(3u^2-4u+4   \right).
\end{align*}
Since this expression is decreasing in $s$, the minimum  is achieved at $s=\frac{1}{2}$. Thus, {\it Case 1} is reduced to {\it Case 2}.

{\it Case 2}.  We now have $s\ge \frac12$.   The left hand side of (\ref{EEq}) decreases when either $\lambda$ or $b$ increase, or when $s$ decreases. Thus, if $s\le 1-\frac{\lambda}{2}$, we can increase $\lambda$, i.e., we can set $1-\frac{\lambda}{2}=s$. If $s\ge 1-\frac{\lambda}{2}$,  we can decrease  $s$, i.e., we can again set $s=1-\frac{\lambda}{2}$. Finally, choosing $b$ as large as possible, we can assume $b=s=1-\frac{\lambda}{2}$. Thus, we  have to check that
\begin{align*}
&\frac{u-2+\frac{1}{u}}{3}-\frac{2s(u^2-u+1)}{3u^2}+\frac{s}{6} +\frac{1}{2}-\frac{u(1-s)}{2}\\
& =
\frac{1}{3u}-\frac{1}{6}-\frac{u}{6}+s\Big(-\frac{2}{3u^2}+  \frac{2}{3u}-\frac{1}{2}+\frac{u}{2}  \Big)
\ge 0.
\end{align*}
Since the coefficient at $s$ is non-negative (recall that $u\ge 1$), it is enough to check the last inequality for $s=\frac{1}{2}$, and it follows from 
$$
-\frac{1}{3u^2}+\frac{2}{3u}-\frac{5}{12} +\frac{u}{12} =\frac{1}{12u^2}(u-1)(u-2)^2\ge 0.
$$

\section{Concluding remarks}

Let $D_n$ be the smallest constant such that for every convex
body $\Omega\subset \mathbb R^n$, every $(n-1)$-dimensional subspace
$H\subset \mathbb R^n$, and every direction $\theta\in H$, we have
$|\langle c(\Omega) - c(P_H \Omega), \theta\rangle|\le D_n
w_\Omega(\theta)$. Here, $P_H\Omega$ is the orthogonal projection of
$\Omega$ to $H$. In \cite{MTY} it was shown that this inequality 
becomes equality for a certain class of convex bodies and an explicit expression was given for $D_n$ as a maximum of a certain rational function of one variable. We suspect that
the same constant $D_n$ should be the supremum of the quantity
$|\langle c(\Omega) - c(\partial \Omega), \theta\rangle|/
w_\Omega(\theta)$ over all convex bodies $\Omega$ and all directions
$\theta$.
However, we do not see how to prove this.

One can also consider other pairs of distinguished points and other notions of ``linear size" of a convex body. In \cite{E}  and \cite{S}, for instance, a sharp upper bound for the distance between the center of mass and the circumcenter  of a  planar convex body
was found in terms of the diameter and the circumradius respectively. However, as far as we know, many natural questions in this area still remain open even in dimension $2$.

{\bf Acknowledgment.} We are grateful to Peter Gordon for the representation of the polynomial at the end of Section \ref{Sec5}.

\end{document}